\documentclass{amsart}
\usepackage[dvips]{graphicx}
\usepackage{amscd}
\usepackage{amsmath}
\usepackage{amsxtra}
\usepackage{amsfonts}
\usepackage{amssymb}

\newtheorem{theorem}{Theorem}[section]
\newtheorem{corollary}[theorem]{Corollary}
\newtheorem{lemma}[theorem]{Lemma}
\newtheorem{proposition}[theorem]{Proposition}

\theoremstyle{definition}
\newtheorem{definition}[theorem]{Definition}
\newtheorem{remark}[theorem]{Remark}

\newtheorem{example}[theorem]{Example}
\theoremstyle{remark}

\renewcommand{\theclaim}{\textup{\theclaim}}

\newtheorem*{acknowledgements}{Acknowledgements}

\numberwithin{equation}{section}

\def\openone%{\hbox{\upshape \small1\kern-3.3pt\normalsize1}}

{\mathchoice

{\hbox{\upshape \small1\kern-3.3pt\normalsize1}}

{\hbox{\upshape \small1\kern-3.3pt\normalsize1}}

{\hbox{\upshape \tiny1\kern-2.3pt\SMALL1}}

{\hbox{\upshape \Tiny1\kern-2pt\tiny1}}}

\makeatletter

\newbox\ipbox

\newcommand{\ip}[2]{\left\langle #1\, , \,#2\right\rangle}
\newcommand{\diracb}[1]{\left\langle #1\mathrel{\mathchoice

{\setbox\ipbox=\hbox{$\displaystyle \left\langle\mathstrut
#1\right.$}

\vrule height\ht\ipbox width0.25pt depth\dp\ipbox}

{\setbox\ipbox=\hbox{$\textstyle \left\langle\mathstrut
#1\right.$}

\vrule height\ht\ipbox width0.25pt depth\dp\ipbox}

{\setbox\ipbox=\hbox{$\scriptstyle \left\langle\mathstrut
#1\right.$}

\vrule height\ht\ipbox width0.25pt depth\dp\ipbox}

{\setbox\ipbox=\hbox{$\scriptscriptstyle \left\langle\mathstrut
#1\right.$}

\vrule height\ht\ipbox width0.25pt depth\dp\ipbox}

}\right. }

\newcommand{\dirack}[1]{\left. \mathrel{\mathchoice

{\setbox\ipbox=\hbox{$\displaystyle \left.\mathstrut
#1\right\rangle$}

\vrule height\ht\ipbox width0.25pt depth\dp\ipbox}

{\setbox\ipbox=\hbox{$\textstyle \left.\mathstrut
#1\right\rangle$}

\vrule height\ht\ipbox width0.25pt depth\dp\ipbox}

{\setbox\ipbox=\hbox{$\scriptstyle \left.\mathstrut
#1\right\rangle$}

\vrule height\ht\ipbox width0.25pt depth\dp\ipbox}

{\setbox\ipbox=\hbox{$\scriptscriptstyle \left.\mathstrut
#1\right\rangle$}

\vrule height\ht\ipbox width0.25pt depth\dp\ipbox}

} #1\right\rangle}

\newcommand{\Tr}{\operatorname*{Tr}}

\newcommand{\cj}[1]{\overline{#1}}

\newcommand{\bz}{\mathbb{Z}}
\newcommand{\M}{\mathcal{M}}

\newcommand{\br}{\mathbb{R}}
\newcommand{\bc}{\mathbb{C}}
\newcommand{\bt}{\mathbb{T}}
\newcommand{\bn}{\mathbb{N}}

\def\blfootnote{\xdef\@thefnmark{}\@footnotetext}

\renewcommand{\mod}{\operatorname{mod}}

\input xy
\xyoption{all}
\usepackage{amssymb}

\def\R{\mathbb{R}}
\def\N{\mathbb{N}}

\def\T{\mathbb{T}}

\def\H{\mathcal{H}}

\def\-{^{-1}}

\def\K{\mathcal{K}}

\def\W{\mathcal{W}}
\def\C{\mathbb{C}}

\begin{document}
\title[Covariant representations]{Covariant representations for matrix-valued transfer operators}
\author{Dorin Ervin Dutkay}
\blfootnote{Research supported in part by a grant from the National Science Foundation DMS 0457491 and the Research Council of Norway, project number NFR 154077/420}
\address{[Dorin Ervin Dutkay]University of Central Florida\\
	Department of Mathematics\\
	4000 Central Florida Blvd.\\
	P.O. Box 161364\\
	Orlando, FL 32816-1364\\
U.S.A.\\} \email{ddutkay@mail.ucf.edu}
\author{Kjetil R\o ysland}
\address{[Kjetil R\o ysland]University of Oslo\\
Department of Mathematics\\
PO Box 1053, Blindern\\
NO-0316 Oslo\\
Norway}

\email{roysland@math.uio.no}
\thanks{
} 
\begin{abstract}
Motivated by the multivariate wavelet theory, and by the spectral theory of transfer operators, we construct an abstract affine structure and a multiresolution associated to a matrix-valued weight. We describe the one-to-one correspondence between the commutant of this structure and the fixed points of the transfer operator. We show how the covariant representation can be realized on $\mathbb{R}^n$ if the weight satisfies some low-pass condition. 
\end{abstract}
\maketitle \tableofcontents
\section{Introduction}\label{intro}

   Several themes are involved in this paper: the ``covariant representations'' in the title has been a central construct from the theory of operator algebras since the 1950s, and they have played a key role in numerous applications since. One of these more recent applications is to a class of wavelets called ``frequency localized'' wavelets. The ``transfer operators'' in the title refers to a construction with origins in probabilistic path models from physics and Ergodic theory.

       One of our aims here is to point out some connections between the two areas, and to show how operator algebraic ideas and representations throw new light on a classical theme.

        Since several ideas are involved, readers from one area may look for pointers to the other.

         We begin with a brief guide to the literature: One use of operator algebras (specifically, $C^*$-algebras) is to the construction of representations. Initially \cite{Mac78}, the focus was on groups, but the notion of covariance from physics (see e.g., \cite{BR87, Rue04}) suggested crossed products of groups with act by automorphisms on $C^*$-algebras (\cite{CMW84, Wil82}). Since the pioneering paper by Stinespring \cite{Sti55}, a preferred approach (e.g., \cite{Arv69}) to constructing representations begins with a positive operator valued mapping, and it was Stinespring's insight that identified the correct ``positivity'' as complete positivity. But, at the same time, notions of positivity are central in a variety of probabilistic path models, beginning with Doeblin \cite{Doe40}, see also \cite{Coh81}. It is now also a key tool in ergodic theory, \cite{Pe83, Wal82}. As a result, Doeblin's operator has subsequently taken on a variety of other incarnations and it is currently known as ``the transfer'' operator, the Ruelle operator, or the Perron-Frobenius-Ruelle operator \cite{Baladi}. The name Ruelle is from its use in statistical mechanics as pioneered by David Ruelle, see \cite{Rue04, Baladi}.

        Of a more recent vintage are applications to wavelets \cite{Dau92}, i.e., special and computational bases in Hilbert space constructed from a class of unitary representation of certain discrete groups of affine transformations. It was realized (e.g., \cite{BraJo02}) that there are transfer operators $R_W$ for wavelets, that the solution to a spectral problem for $R_W$ yields wavelet representations; and moreover that these representations come with a useful covariance. Hence the circle closes with the positivity question from \cite{Sti55}, albeit in a different guise. It is intriguing that a variant of the operators $R_W$ have now also found use in quantum error correction codes, see \cite{CKZ06, Kri06}.

        All the versions of transfer operators involve hierarchical processes with branching, and probabilities assigned by a weight function. Our paper focuses on $R_W$ in the wavelet context; and we demonstrate that for many wavelets, the weight function must take the form of an operator transformation  $X \mapsto m X m^*$ where $m$ is a fixed matrix valued function and where $m^*$ denotes the adjoint operator, in this case transpose-conjugate. This matrix version of $R_W$ is necessary for understanding wavelet constructions associated to wavelet sets \cite{DaLa98}, and more generally to non-MRA wavelets, \cite{BCM02}. (MRA stands for multiresolution analysis \cite{Dau92}.)

 An orthogonal {\it wavelet} is a function $\psi\in L^2(\br)$ such that $\{2^{j/2}\psi(2^j\cdot-k)\,|\,j,k\in\bz\}$ is an orthonormal basis in $L^2(\br)$. The main technique used to construct wavelets is by {\it multiresolutions}. For this one needs a {\it low-pass filter} $m_0$, i.e., a $2\pi$-periodic Lipschitz continuous function that satisfies the low-pass condition and the QMF condition:
 $$m_0(0)=\sqrt{2},$$
 $$\frac{1}{2}\left(\left|m_0\left(\frac{x}{2}\right)\right|^2+\left|m_0\left(\frac{x+2\pi}{2}\right)\right|^2\right)=1,\quad(x\in\br).$$
  Then, from the low-pass filter, one constructs the Fourier transform of the {\it scaling function} 
  $$\hat\varphi(x)=\prod_{n=1}^\infty\frac{m_0\left(\frac{x}{2^n}\right)}{\sqrt{2}},\quad(x\in\br).$$
  Finally, the wavelet is constructed from the scaling function using the formula
  $$\hat\psi(x)=\frac1{\sqrt{2}}e^{i\frac x2}\cj m_0\left(\frac{x}{2}+\pi\right)\hat\varphi\left(\frac{x}{2}\right),\quad(x\in\br).$$
  
  It turns out that the low-pass condition and the QMF condition, while necessary, are not always sufficient to obtain an {\it orthonormal} wavelet. One of the extra conditions on $m_0$ that guarantees the orthogonality of the wavelet was given by W. Lawton \cite{Law91}. Here the {\it transfer operator} was introduced in the study of wavelets. 
  
  In this context, the transfer operator is defined on $2\pi$-periodic functions by:
  $$R_{m_0}f(x)=\frac{1}{2}(|m_0(\frac x2)|^2f(\frac x2)+|m_0(\frac{x+2\pi}{2})|^2f(\frac{x+2\pi}{2})),\quad(x\in\br).$$

 Lawton's condition states that the wavelet is orthogonal if and only if the only continuous functions $h$ with $R_{m_0}h=h$ are the constants. 
 
 When this condition is not satisfied, the resulting wavelet still has an interesting property, namely it generates a Parseval frame. Having this, the theory of dilations of Parseval frames due to D. Han and D. Larson \cite{HaLa00} can be used. The wavelet Parseval frame is the projection of an orthonormal basis in a bigger space. The problem then was if this orthonormal basis has a similar wavelet structure, i.e., if it is generated by the application of unitary dilation and translation operators $U$ and $T$ that satisfy the commutation relation $UTU^{-1}=T^2$ (resembling the relation between the dilation and the translation operators in $L^2(\br)$). The answer to this question is positive and it is given in \cite{BDP}. It turns out that each fixed point of the transfer operator $R_{m_0}$ will give rise to such a wavelet structure, and we call this a {\it covariant representation}. Putting together these representations, one obtains the covariant representation which has an orthonormal wavelet in a bigger Hilbert space whose projection onto $L^2(\br)$ is the Parseval wavelet frame constructed in the classical wavelet theory. This is one of the wonderful uses of covariant representations. Since these orthonormal wavelet bases live in a bigger Hilbert space, i.e., one that contains $L^2(\br)$ as a subspace, Han and Larson coined the term ``super-wavelets''. 
  Another interesting application of covariant representation is the computation of the peripheral spectrum for the transfer operator: the fixed points of $R_{m_0}$ are in one-to-one correspondence with the commutant of the covariant representation. Often this commutant can be explicitly computed, hence the eigenspace of $R_{m_0}$ can be obtained from that. The spectral properties of transfer operators play an important role in the ergodic analysis of discrete dynamical system (see \cite{Baladi}).
  
  Why matrix-valued transfer operators? It is known (see \cite{Dau92}) that not all wavelets in $L^2(\br)$ can be constructed from a multiresolution. However there are generalizations of this that will do the job: each orthogonal wavelet can be constructed from a {\it generalized multiresolution analysis}, a notion introduced by Baggett et al. \cite{BMM99}. The construction requires some matricial low-pass filters. Also multiwavelet theory, in which one uses more than one function to generate the basis, requires matricial filters. But many times the resulting wavelet is only a Parseval frame, not an orthogonal basis (see \cite{BJMP05}).
  
  There is no analogue for the Lawton condition in the case of matricial low-pass filters. The main reason for this is the impossibility to give a good generalization of the notion of zeros for the fixed points of the transfer operator, a notion which plays an essential role in the scalar case. We believe that the covariant representations can provide a way around this, and we can analyze the orthogonality of wavelets and scaling functions constructed from matricial low-pass filters through a study of the associated covariant representations.

%%%%%%%%%%%%%%%%%%%%%%%%%%

More generally, a transfer operator, also called Ruelle operator, is associated to a finite-to-one continuous endomorphism on a compact metric space $r:X\rightarrow X$ and a  weight function $W:X\rightarrow[0,\infty)$, and it is defined by
$$R_Wf(x)=\sum_{r(y)=x}W(y)f(y)$$
for functions $f$ on $X$. 
\par
Transfer operators have been extensively used in the analysis of discrete dynamical systems \cite{Baladi} and in wavelet theory \cite{BraJo02}.
\par
In multivariate wavelet theory (see for example \cite{Jiang} for details) one has an expansive $n\times n$ integer matrix $A$, i.e., all eigenvalues $\lambda$ have $|\lambda|>1$, and a multiresolution structure on $L^2(\br^n)$, i.e., a sequence of subspaces $\{V_j\}_{j\in\bz}$ of $L^2(\br^n)$ such that
\begin{enumerate}
\item $V_j\subset V_{j+1}$, for all $j$;
\item $\cup_j V_j$ is dense in $L^2(\br^n)$;
\item $\cap_j V_j=\{0\}$;
\item $f\in V_j$ if and only if $f( (A^T)^{-1}\cdot)\in V_{j-1}$;
\item There exist $\varphi_1,...,\varphi_d\in V_0$ such that $\{\varphi_k(\cdot-j)\,|\,k\in\{1,...,d\},j\in\bz^n\}$ forms an orthonormal basis for $V_0$.
\end{enumerate}

The functions $\varphi_1,...,\varphi_d$ are called {\it scaling functions}, and their Fourier transforms satisfy the following {\it scaling equation}:
$$\hat\varphi_i(x)=\sum_{j=1}^dq^{-1/2}m_{ji}(A^{-1}x)\hat\varphi_j(A^{-1}x),\quad(x\in\br^n,i\in\{1,...,d\}),$$
where $q:=|\det A|$ and $m_{ji}$ are some $\bz^n$-periodic functions on $\br^n$. 
\par
The orthogonality of the translates of $\varphi_i$ implies the following {\it QMF equation}:
$$\frac1q\sum_{Ay=x\mod\bz^d}m^*(y)m(y)=1,\quad(x\in\br^n/\bz^n),$$
where $m$ is the $d\times d$ matrix $(m_{ij})_{i,j=1}^d$. 
When the translates of the scaling functions are not necessarily orthogonal, one still obtains the following relation: if we denote by
$$h_{ij}(x):=\sum_{k\in\bz^n}\cj{\hat\varphi_i}(x+k)\hat\varphi_j(x+k),\quad(x\in\br^n/\bz^n),$$
then the matrix $h=(h_{ij})_{i,j=1}^d$ satisfies the following property:
\begin{equation}\label{eqrh}
Rh(x):=\frac1q\sum_{Ay=x\mod\bz^n}m^*(y)h(y)m(y)=h(x),\quad(x\in\br^n/\bz^n),
\end{equation}
i.e., $h$ is a fixed point for the {\it matrix-valued transfer operator
}$R$. The fixed points of a transfer operator are also called {\it harmonic functions} for this operator. Thus the orthogonality properties of the scaling functions are directly related to the spectral properties of the transfer operator $R$. 
\par
This motivates our study of the harmonic functions for a matrix-valued transfer operator. The one-dimensional case (numbers instead of matrices) was studied in \cite{BraJo02,Dugale,Duspec}. These results were then extended in \cite{DuJoifs,DuJomart}, by replacing the map $x\mapsto Ax\mod\bz^n$ on the torus $\bt^n/\bz^n$, by some expansive endomorphism $r$ on a metric space.
\par
Here we are interested in the case when the weights defining the transfer operator are matrices, just as in equation \eqref{eqrh}. We keep a higher level of generality because of possible applications outside wavelet theory, in areas such as dynamical systems or fractals (see \cite{DuJoifs,DuJomart}). However, for clarity, the reader should always have the main example in mind, where $r:x\mapsto Ax\mod\bz^n$ on the torus $\bt^n$.  
\par
In \cite{DuJomart} it was shown that, whenever a pair $(m,h)$ is given, with $h\geq 0$ and $Rh=h$, one can construct a covariant structure on some Hilbert space $H$, i.e., a unitary $U$, a representation $\pi$ of continuous functions on $X$, and some scaling functions $\varphi_1,\dots,\varphi_d\in H$ such that 
$$U\pi(f)U^*=\pi(f\circ r),\quad(f\in C(X)),$$
$$\ip{\varphi_i}{\pi(f)\varphi_j}=\int_Xfh_{ij}\,d\mu,\quad(f\in C(X)),$$
\begin{equation}\label{eqscal}U\varphi_i=\sum_{j=1}^d\pi(m_{ji})\varphi_j,\quad(i\in\{1,\dots,d\}),
\end{equation}
\begin{equation}\label{eqcorr}
\{U^{-n}\pi(f)\varphi_i\,|\,n\in\bz,i\in\{1,\dots,d\},f\in C(X)\}\mbox{ is dense in }H,
\end{equation}
where $\mu$ is a strongly invariant measure on $X$ (see \eqref{strong invariance}).
\par
Moreover, this construction is unique up to isomorphism.
Thus, from the ``filter'' $m$ and the ``harmonic function'' $h$, one can construct a multiresolution structure, similar to the one used in wavelet theory (see \cite{BraJo02}). There are several uses for such a structure: one can analyze the peripheral spectrum of the transfer operator \cite{Dugale}, construct super-wavelet bases on spaces bigger than $L^2(\br)$ or on some fractal spaces \cite{BDP,DuJowf}, or use the rich algebraic and analytic structure of this multiresolution to perform computations needed in the harmonic analysis of fractal measures \cite{DuJoifs}. 
\par
This covariant representation was analyzed in more detail in the one-dimensional case $d=1$ in \cite{DuJomart,DuJoifs}. The main tool used there was the introduction of some random-walk measures $P_x$ following an idea of Conze and Raugi \cite{MR1078079}. 
\par
In our case $m$ and $h$ are matrix-valued. We will use the language of vector bundles because we are interested also in continuous harmonic functions and projective multiresolution analyses (see \cite{DuRo1,PackerRieffel1}). We will describe the covariant representation in the multivariable case using some positive-matrix valued measures (see Theorem \ref{representation}, Propositions \ref{prop4-3} and \ref{proppx}), thus extending previous results from \cite{DuJomart}. 
\par
In Section \ref{def} we introduce the setup and the main definitions and assumptions on the transfer operator $R$ associated to a matrix-valued weight $h$. Then in Section \ref{cova} we show how the covariant representation can be constructed on the solenoid of $r$. In Section \ref{harm} we show that the operators that commute with this covariant representation are in one-to-one correspondence with the harmonic functions of $R$ (Theorem \ref{cone}). Then the set of harmonic functions inherits a $C^*$-algebra structure from the commutant, and a more intrinsic description of the multiplication is given in Theorem \ref{thprod} and Remark \ref{remprodcont}.\par
The covariant structure can be described in terms of some operator-valued measures $P_x$. We do this in Section \ref{cocy}. With the aid of these measures we can give another form to the correspondence between harmonic maps and operators in the commutant, which can be identified now with {\it cocycles} (Proposition \ref{conv}, Theorem \ref{expe}). 

\par
In the case when the filter $m$ satisfies a low-pass condition (in this case, the $E(l)$-condition as defined in \cite{Jiang}), we can realize the covariant representation in the more familiar environment of $L^2(\br^n)$. The space $\br^n$ has a natural embedding in the solenoid, the measures $P_x$ are atomic and, the measure of the atoms are directly related to the solutions of the refinement equation, i.e., the scaling functions. Our results extend some ideas from \cite{MR1777126,MR1794571,MR2254502} to the matricial case.

\section{Definitions and preliminaries} \label{def}
\par
{\bf The dynamical system.}
Let $X$ be a compact Hausdorff space with a surjective  and finite to one
continuous map $r: X \rightarrow X$. Moreover, let $\mu$ be a regular 
measure on $\mu$ that is strongly $r$-invariant, i.e. 
\begin{equation}
  \label{strong invariance}
  \int_X  f d \mu =   \int_X  \frac1{\# r\-(t)} \sum_{rs=t} f(s) d\mu(t),
\end{equation}
for every $f \in C(X)$.\par
{\bf Main example.}
The main example we have in mind,  is the following: 
Let $G \subset \R^n$ be a discrete subgroup  such that  $\R^n /G$ is compact, i.e.,
$G$ is a full-rank lattice and $\R^n /G \simeq \T^n=:X$.
Moreover, let $A \in GL(\R^n)$  be strictly expansive and such that 
$A G \subset G$.
Let $p : \R^n \rightarrow X$ denote the quotient map and define  
a map $r : X \rightarrow X$ as 
$r(p(x)) = p(Ax)$. This is a  $|\det(A)|$-folded normal covering  map. 
Finally, we let  $\mu$ be the Haar measure on $\T^n$. 

\par
{\bf The transfer operator.}
Let $\rho : \xi \rightarrow X$ be  a $d$-dimensional complex vector bundle over $X$ with
a Hermitian metric on $\xi$, i.e., a continuous map 
$\langle \cdot , \cdot \rangle : \xi \times \xi \rightarrow \C$ that restricts to 
an inner-product on each fiber. Such a map always exist when $X$ is compact, 
\cite[1.3.1]{Atiyah1}. Let $S$ be the set of continuous sections in $\xi$. Then
$$
s_1, s_2 \mapsto \int_X \langle s_1(x) , s_2(x) \rangle d\mu(x),
$$
defines an inner-product on $S$. Let $\K$ denote the Hilbert space obtained by the completion with respect to the corresponding norm. Then $\K$ 
is the Hilbert space  of  $L^2$-sections in $\xi$ with respect to 
the measure $\mu$.
\par
The $C(X)$-product on the sections in $\xi$ gives a representation of $C(X)$ on $\K$ by pointwise multiplication,
$\kappa : C(X) \rightarrow B(\K)$ such that 
$(\kappa(f)s)(x) = f(x) s(x)$ for every $x \in X$. Every  bundlemap on 
$\xi$ commutes with this representation.
\begin{remark}\label{remcomm}
Let us consider the von Neumann algebra $\kappa(C(X))''$ generated by this representation.
Since the bundle is locally trivial around any point in $X$, we can take $U_1, \dots, U_r$ open sets in $X$,
with bundle isomorphisms $\phi_j : \xi|_{U_j} \rightarrow U_j \times \C^d$. 
\par
If $V \in \kappa(C(X))''$ there exists $V(x) \in \C$ such that 
$x \mapsto \phi_j(x)\- V(x) \phi_j(x)$ is 
 an $L^\infty(U_j, \mu|_{U_j})$ function
 and $\phi_j(Vs (x)) = \phi_j(V(x)s(x))$, $\mu$-a.e.
\par
If $V \in \kappa(C(X))'$ there exists $V(x) \in \text{End}_{\C} (\xi|_x)$ 
such that the map $x \mapsto   \phi_j(x)\-V(x)\phi_j(x)$ is contained in 
$ M_d(\bc)\otimes L^\infty(U_j,\nu|_{U_j})$ and 
$\phi_j(Vs) (x)) = \phi_j(V(x)s(x))$ $\mu$-a.e.
\end{remark}
\par
Let $r^*\xi$ denote the pull-back of $\xi$ along $r$ (see \cite[1.1]{Atiyah1}).
The space of sections in $r^*\xi$  is endowed with the pull-back Hermitian metric from $\xi$. 
Let $\tilde \K$ denote the $L^2$-sections in this bundle with respect to this metric 
and $\mu$. Moreover, let $\tilde \kappa: C(X) \rightarrow B(\tilde \K)$ 
denote the corresponding representation of $C(X)$ by pointwise multiplication.  
\par
The {\it weight} (or the {\it filter} if we use wavelet terminology) that is used to define the transfer operator is in our case an operator $m \in B( \tilde \K , \K)$
such that $m \tilde\kappa(f) = \kappa(f) m $ for every $f \in C(X)$. As in Remark \ref{remcomm}, $m$ is a pointwise multiplication by a linear map between the fibers $\xi_{rx}$ and $\xi_x$, i.e., there exist 
a unique $m(x) \in \text{Hom}_\C(\xi|_{rx}, \xi|_{x})$
such that $ \phi_j(m(x,v)) = \phi_j(m(x)v)$ for every $x \in U_j$ and $v \in \xi_{rx}$.
\begin{definition}
Let $\M = \kappa(C(X))'$. As we have seen before, the space $\M$ consists in bounded measurable bundle maps on $\xi$. We define the {\it transfer operator} on $\M$ by
$$
(Rf) (x) = \frac 1 {\#r^{-1}(x)} \sum_{ry = x}    m^* (y) f (y)   m (y),\quad(x
\in X).
$$ 
\end{definition}
\par
{\bf Assumptions on the filter $m$.}
Throughout the paper we will assume that $m$ satisfies the following conditions:
\begin{enumerate}
  \item $m \in B(\tilde \K , \K)$ is injective\label{injec}
  \item $\sup_k  \|R^k \| < \infty$ \label{bound }
  \item There exists an $h \in \text{End}(\xi)$ such that 
   $Rh = h$ and $h\geq 0$ (here we mean $h$ non-negative as an operator). 
\end{enumerate} 
\begin{remark} The first condition \eqref{injec} implies  that $m(x)$ is
invertible $\mu$-a.e.
\par
Note that if there exists an  $h\in\M$ such that $Rh=h$ and $h\geq c1$ for
some $c>0$ (here $1$ is the identity bundle map on $\xi$), then
\eqref{bound } holds automatically, see \cite{DuRo1}.
\end{remark}
\begin{definition}
A function $f\in\M$ such that $Rf=f$ is called {\it harmonic} with respect to the transfer operator $R$. We denote by $\mathfrak H$ the set of all bounded harmonic functions. 
\end{definition}

\section{The covariant representation associated to a matricial weight and a harmonic bundle map}\label{cova}
 We mentioned in the introduction that for some choices of the filter $m_0$ such as the ``stretched Haar filter'' (see \cite{Dau92}) the scaling function and wavelet constructed in $L^2(\br)$ are not orthogonal,i.e., they do not have orthogonal translates. Thus classical wavelet theory breaks down if we want to have orthogonal solutions for the prescribed scaling equation. Even deeper problems appear when the filter does not satisfy the low-pass condition, since the infinite product defining the scaling function from $m_0$ is $0$. Fractal spaces may occur (see \cite{DuJowf}). We will show that we always have a solution to a given scaling equation in some Hilbert space with a covariant representation. The Hilbert space depends sensitively on the filter $m$ and the correlation function $h$.  
 
 We construct now the affine structure  associated to the matricial filter $m$ and the positive harmonic function $h$. The covariant representation is a Hilbert space endowed with an affine structure given by a unitary $U$, which takes the place of the dilation operator, and a representation $\pi$ of $C(X)$, which takes the place of the representation generated by translations. In this Hilbert space one has several vectors $\varphi_1,\dots,\varphi_n$ that satisfy the scaling equation with the prescribed filter $m$ as in \eqref{eqscal}, and having correlation function $h$ as in \eqref{eqcorr}.

{\bf The ground space.}
Consider the projective system
$$\label{projective system}
\xymatrix{
X & X \ar[l]_r & X \ar[l]_r& \cdots \ar[l] \\ 
}
$$  
and let $X_\infty$ be its projective limit, 
with  projections $\theta_n : X_\infty \rightarrow X$, $(n\geq0)$, and a homeomorphism 
$\hat r: X_\infty \rightarrow X_\infty$ such that the following diagrams are commutative for all $n\geq0$:
$$
\xymatrix{
X_\infty \ar[d]_{\theta_n} \ar[r]^{\hat r} \ar[dr]^{\theta_{n-1}} 
& X_\infty \ar[d]^{\theta_n} \\
X \ar[r]_r & X \\
}
$$
We have the following identification of $X_\infty$
$$
X_\infty = \{( x_0 , x_1 , x_2, \dots) \in X \times X \times \dots | rx_{j+1} =  x_j,\mbox{ for all }j\geq0     \}
$$
$$\theta_k (x_0, x_1, \dots ) = x_k,\quad \hat r(x_0,x_1,\dots)=(rx_0,x_0,x_1,\dots).$$  
The space $X_\infty$ is a compact Hausdorff space with the topology generated by 
the inverse images of the open sets in $X$ with respect to the maps $\theta_n$, $n \geq 0$.
\par
{\bf The Hilbert space.}
Let $S_k$ denote the vector space of continuous sections in 
$\theta_0^*\xi$ that depend only on the first $k+1$ coordinates $x_0,x_1,\dots,x_k$. We see that 
$S_k$ is formed by $s \circ \theta_k$ for sections in $\xi$, $s: X \rightarrow \xi$, i.e.,
maps of the form 
$$x \mapsto (\hat r ^{-k} (\tilde x) , s \circ \theta_k(\tilde x) ),\quad(\tilde x\in X).$$

We define the inner-product $\langle \cdot , \cdot \rangle_k : S_k \times S_k \rightarrow \C$
by
$$ \langle f \circ \theta_k , g\circ \theta_k  \rangle_k  
  :=  \int_{X}   \frac 1 {\#r^{-k}(x)} \sum_{r^ky = x} 
 \langle m^{(k)}(y)  
f(y ) ,  h (y   ) m^{(k)}(y)  
g(y ) \rangle\,d\mu(x),
$$
where 
$$m^{(k)}(x) := m ( x) m(rx)  \dots    
 m(r^{k-1}x).$$
We have the following compatibility relation between these inner-products: 
    \begin{align*}
      \langle f \circ \theta_l , g \circ \theta_l  \rangle_{k+l}  =& \ip{f\circ r^k\circ\theta_{k+l}}{g\circ r^k\circ\theta_{k+l}}_{k+l}\\
      =&  
      \int_X  \frac 1 {\#r^{-(k+l)}(x)}       
       \sum_{r^{k+l} y = x} 
     \langle m^{(k+l)}(y)  
f(r^ky ) ,   h ( y  )  m^{(k+l)}(y)  
g(r^ky) \rangle d\mu(x)
      \\
   = &  \int_X 
   \frac 1 {\#r^{-l}(x)}  \sum_{r^l y= x} \langle m^{(l)}(y)  
f(y) ,  R^k(h) ( y  )
m^{(l)}(y)  
g(y) \rangle d\mu(x)     \\
= & \langle f \circ \theta_l  , g \circ \theta_l \rangle_l\quad(\mbox{since }R^kh=h).
    \end{align*}
This shows that the restriction $\langle \cdot , \cdot \rangle_{k+l} | _{S_l \times S_l}  = \langle \cdot , \cdot \rangle_l$.

We obtain a sesquilinear form $\ip{\cdot}{\cdot}$ on 
$\cup_k S_k$, i.e., a (possibly degenerate) inner-product. Let $\H$ denote the 
Hilbert space completion with respect to this inner-product. 
\par
An application of the Stone-Weierstrass theorem gives that 
$\cup_k S_k$ is dense in the space of continuous sections in $\theta_0^* \xi$. Moreover, the continuous sections in $\theta_0^*\xi$ are contained in $\H$ by the following argument.  
\begin{lemma}
 If $\{s_k\}_{k \in \N}$ is
a sequence of sections in $\xi$ such that $\{s_k\circ\theta_k\}_{k\geq0}$
converges uniformly to a section $s : X_\infty \rightarrow \theta^*_0 \xi$
then  $\langle s, s \rangle = \lim_k \langle s_k\circ\theta_k , s_k \circ
\theta_k \rangle_k$.
\end{lemma}

\begin{proof}
 \begin{align*}
 &\| s_{k+l} \circ \theta_{k+l}  - s_k \circ \theta_k \|^2 \\
%  & \| (s_{k+l}   - s_k \circ r^l )\circ  \theta_k \|^2 \\
  =&     \int_X  \frac 1 {\#r^{-(k+l)}(x)}
      \sum_{r^{k+l} y = x}
    \langle
s_{k+l} (y) - s_k \circ r^l(y) ,   m^{(k+l)*}(y)  h ( y  )  m^{(k+l)}(y)
(s_{k+l} (y) - s_k \circ r^l(y)) \rangle d\mu(x)\\
       \leq&   \| s_{k+l}\circ  \theta_{k+l} - s_k \circ \theta_k
\|^2_\infty
        \int_X  \frac 1 {\#r^{-(k+l)}(x)}  \sum_{r^{k+l} y = x} \|
m^{(k+l)*}(y)  h ( y  )  m^{(k+l)}(y)\| d \mu(x)\\
      \leq&   \| s_{k+l}\circ  \theta_{k+l} - s_k \circ \theta_k
\|^2_\infty   d \|R^{k+l} h \|_\infty,
      \end{align*}
 and the last inequality follows since
  $\sum_i \| H^*_i H_i \| \leq \sum_i \text{Tr}(H^*_i H_i) = \text{Tr} (
\sum_i  H_i^* H_i)
 \leq d \|\sum_i H_i^* H_i \|  $ for every finite sequence of matrices
$H_i$.
\end{proof}

Next we define the covariant representation and the multiresolution on the Hilbert space $\H$.

\begin{theorem} \label{representation}
Let  $P_k$ denote the orthogonal projection in $B(\H)$ onto the closed subspace $\H_k$ generated by $S_k$. Define $$\pi : C(X) \rightarrow B(\H),\quad 
(\pi(f)s)(x) = f \circ \theta_0(x) s(x),\quad(f\in C(X),s\in\H,x\in X_\infty);$$ 
$$U \in B(\H),\quad(U s )(x) = m \circ \theta_0(x) s \circ \hat r(x),\quad(s\in\H,x\in\tilde X).$$ The following relations hold: 
\begin{enumerate}
  \item \label{unitary} 
    $U$ is a unitary such that $U \pi(f) U^* = \pi(f \circ r)$ for all $f\in C(X)$.
  \item \label{commutant}
    $P_0 \in \pi(C(X))'$.
  \item \label{increasing}
    $P_k \leq P_{k+1}$ for every $k \in \N$.
  \item \label{dense}
    $\lim_k  P_k s  =  s $ for every $s \in \H$. 
  \item \label{relation} 
    $UP_{k+1} U^* = P_k$ for every $k \in \N$.
\end{enumerate}
\end{theorem}

\begin{proof}
  \eqref{commutant} and \eqref{increasing} follow directly from the definition. 
\eqref{dense}  follows from the fact that $\cup_k  S_k$ is dense in
the sections  of $\theta_0^* \xi$.
  
We must  prove \eqref{unitary}. 
First, let $f, g \in S$, then 
\begin{align*}
  & \langle m \circ \theta_0 f\circ \theta_k\circ \hat r , m\circ\theta_0g\circ \theta_k\circ\hat r  \rangle = 
   \langle (m \circ r^{k-1}  f) \circ \theta_{k-1}, (m \circ r^{k-1}  g) \circ \theta_{k-1}
   \rangle\\
        = & \int_X \frac 1 {\# r^{-(k-1)}(x)} \sum_{r^{k-1} y = x}  
   \langle m^{(k-1)} ( y ) m(r^{k-1}y) f (y),  h ( y   )   
   m^{(k-1)} ( y) m(r^{k-1}y) g ( y) \rangle d \mu(x)  
\\ 
 = & \int_X  \frac 1 {\# r^{-(k-1)}(x)}   \sum_{r^k x = y } 
    \langle m^{(k)} ( y)  f (y), h (y   )  m^{(k)} (y)g ( y) \rangle d \mu(x)
    \\ 
    = &  \langle f \circ \theta_k , g \circ \theta_k \rangle
\end{align*}
This shows that $U$ is isometric  on  $\cup_k S_k$.

To see that $U$ is surjective on $\H$, we will show that $U\H_{k+1} = \H_k$, (with 
$\H_k:=P_k\H$) which by \eqref{dense} implies that $U$ is surjective.  
If $s \in \H_{k+1}$, then $U s = m\circ \theta_0  s \circ \hat r \in\H_k$. 
If we can show that $\{ ( m \circ r^k s) \circ \theta_k | s \mbox{ measurable section in }\xi\}$ 
is dense in  $\H_k$, we are done. Recall that $m(x)$ is   invertible $\mu$-a.e. 
\par
Take $s\in S$. Suppose, $\{V_l\}_{l \in \N}$ is a decreasing family of open sets such that $ \{ x \in X |  \det m\circ r^k (x) = 0\} \subset V_l$  
for every $l \in \N$ and $\lim_l \mu (V_l ) = 0$. (Note also that by \eqref{strong invariance} the measure $\mu$ is invariant under $r$, i.e., $\mu(r^{-1}(E))=\mu(E)$ for all measurable sets $E$). For every $l \in \N$, there 
exists  a measurable  section $s_l$ in $\xi$,  such that 
$s_l$ is $0$ on $V_l$, and $m\circ r^k (x) s_l (x) = s (x)$ for $x\in X\setminus V_l$.  
Then $\lim_l m\circ r^k s_l =  s$, $\mu$-a.e. and by the dominated convergence theorem, $\lim_l Us_l \circ \theta_{k+1}  =  s \circ  \theta_k$, in $\H$, i.e.,
$s \circ \theta_k$ is in the closure  of $U S_k$.  
 \par
 Finally \eqref{relation} follows from the relation  $U\H_{k+1} = \H_k$ and \eqref{unitary}.
\end{proof}

\begin{definition}  We denote by $C^*(X,r,m,h)$ the $C^*$-algebra generated by $U$ and $\pi(f)$, $f \in C(X)$, and we call it the covariant representation associated to $m$ and $h$.
\end{definition}

\section{Harmonic and measurable bundle maps }\label{harm}
In this section we establish the one-to-one correspondence between harmonic functions and operators in the commutant $C^*(X,r,m,h)'$. We begin with a result on self-adjoint elements, and we require that the harmonic function $h_0$ be dominated by $h$, in the sense that $|h_0|\leq ch$ for some constant $c>0$. In the case when $h$ is bounded away from $0$, i.e., $h\geq C1$ for some $C>0$, then this domination condition is automatically satisfied and we can extend the correspondence to non-self-adjoint elements (Corollary \ref{corcore}).

\begin{theorem} \label{cone}
{\rm (i)}{\bf [Harmonic maps to operators in the commutant]}
  If $h_0 \in \M$, $h_0^* = h_0$,  $Rh_0 = h_0$ and there exists a positive number $c \geq 0$ 
  such that $|h_0| \leq c h$, there exists a unique selfadjoint operator $A$ in the commutant $C^*(X,r,m,h)'$ such that 
  $$\ip{f\circ\theta_0}{Ag\circ\theta_0}=\int_X\ip{f}{h_0g},\quad(f,g\in S).$$ We denote this operator by $A_{h_0}$.
  Moreover in this case
  $$
  \langle f\circ \theta_k, A g \circ \theta_k \rangle = 
  \int_X \frac{1} {\# r^{-k} (x)} \sum_{r^k y = x}  \langle m^{(k)}(y) f (y), h_0(y) m^{(k)}(y) g(y) \rangle \mu(dx),
  $$
  for every $k \in \N$ and $f,g\in S$. 
  \par
  {\rm (ii)}{\bf [Operators in the commutant to harmonic maps]}
 Conversely, define the operator $T \in B( P_0 \H, \K)$ such that $T s \circ \theta_0 = h s$ for every $s \in S$.
If $A$ is a self-adjoint operator in the commutant $C^*(X,r,m,h)$ then $h_A:=TP_0AP_0T^*\in\M$ is a harmonic function, i.e., $Rh_A=h_A$, and $|h_A|\leq ch$ for some constant $c>0$. Moreover, the correspondences described in (i) and (ii) are inverses to each other, i.e.,
$$h_{A_{h_0}}=h_0,\quad A_{h_A}=A.$$
  \end{theorem}

  \begin{proof}
  (i)  Let $B_k$ denote the sesquilinear form on $S_k$, given by
    $$ B_k( f \circ \theta_k , g \circ \theta_k )   = 
    \int_X \frac{1} {\# r^{-k} (x)} \sum_{ r^k y = x}  \langle m^{(k)}(y) f (y), h_0(y) m^{(k)}(y) g(y) \rangle \mu(dx).
    $$
    Using the  computation for $\langle \cdot, \cdot \rangle_k$, 
    we see that $B_k | _{S_{k-1} \times S_{k-1}} = B_{k-1}$.
    \par
    By the boundedness assumption on $h_0$, we have $-ch\leq h_0\leq ch$. This shows that
    $$ -c\ip{f\circ\theta_k}{f\circ\theta_k}\leq B_k(f\circ\theta_k,f\circ\theta_k)\leq c\ip{f\circ\theta_k}{f\circ\theta_k}.$$
    Since $h_0$ is self-adjoint this implies that we obtain a bounded sesquilinear map $B$ on $\H$
    that restricts to $B_k$ on $S_k$.  Let $A$ denote the bounded operator on $\H$ 
    such that $B(f,g) = \langle f,Ag \rangle$ for every $f,g \in \H$. 
  
    The computation that showed that $U$ was isometric with respect to $\langle \cdot, \cdot \rangle$ applies here too, and it shows that $B(f,g) = B(Uf, Ug)$ for every $f,g \in \H$. Now 
    $$
    \langle f, U A g \rangle = \langle U^* f, Ag \rangle = B( U^*f , g) = B(f, U g) =  \langle f, A U g \rangle,
    $$
    i.e., $A U = UA$. Moreover, by a direct computation, we have $B(f, \pi(a) g ) = B( \pi(a^*) f, g)$ 
    for every $a\in C(X)$, $f,g \in \H$. Using this, we obtain 
    $$
    \langle f, A \pi(a) g \rangle = B( f, \pi(a) g) = B(\pi(a^*)f, g) = \langle \pi(a^*) f, A g \rangle,
    $$
    i.e.,  $\pi(a) A = A\pi(a)$  for every $a \in C(X)$.
    Finally, since $U$ is unitary, $U^* A = U\- A = A U\- = AU^*$ and $A$ commutes with every 
    operator in $C^*(X,r,m,h)$.
    \par
    The uniqueness follows from the fact that $\cup_kS_k$ is dense in $\H$.
    
\par   
(ii) 
 First, a simple computation shows that $T^*s=s\circ\theta_0$ for $s\in\K$. Then
$$
\int_X \langle s_1(x) , (T a T^* )(x)s_2(x) \rangle d\mu(x)   
= \langle s_1 \circ \theta_0 , a s_2 \circ \theta_0    \rangle.
$$
for every $a \in \pi(C(X))' \cap B(P_0 \H )$. 

If $A$ commutes with $\pi(f)$ for all $f\in C(X)$ then, since $P_0$ commutes also with $\pi(f)$, it follows that $P_0AP_0$ commutes with $\pi(f)$ so $TP_0 A P_0 T^*$ commutes with $\kappa(f)$ for all $f\in C(X)$. Therefore $TP_0AP_0T^*\in  \kappa(C(X))'$, and we can use Remark \ref{remcomm}. Moreover, if
$A \in C^*(X,r,m,h)'$, it commutes also with $U$, and then, using the strong invariance of the measure $\mu$, 
\begin{align*}
& \int_X \langle s_1 , R( T P_0 A P_0 T^*) s_2 \rangle d \mu = 
\int_X  \frac{1}{\#r^{-1}x} \sum_{ry = x} \langle m(y) s_1(ry) ,  (T P_0 A P_0 T^*)(y)  m(y) s_2(ry) \rangle d \mu(x) \\
= & \int_X   \langle m (x)  s_1(rx) ,
(T  P_0 A P_0T^* )(x)  m(x) s_2(rx) \rangle d \mu(x) 
=  \langle Us_1\circ \theta_0  , P_0 A P_0 Us_2 \circ \theta_0 \rangle\\ 
= & \langle Us_1 \circ \theta_0 , A Us_2 \circ \theta_0  \rangle 
=  \langle s_1 \circ \theta_0 , A s_2 \circ \theta_0 \rangle 
=  \langle s_1\circ \theta_0 ,  P_0 A P_0 s_2\circ \theta_0  \rangle \\
= & \int_X \langle s_1 , T P_0 A P_0 T^* s_2 \rangle d \mu,
\end{align*} 
i.e., $T P_0 AP_0 T^*$  is a  fixed point for $R$.
\par
Since $A$ is bounded and self-adjoint, it follows that 
$$|\int_X\ip{s_1(x)}{h_A(x)s_1(x)}\,d\mu(x)|=|\ip{s_1\circ\theta_0}{As_1\circ\theta_0}|$$$$\leq\|A\|\ip{s_1\circ\theta_0}{s_1\circ\theta_0}=\|A\|\int_X\ip{s_1(x)}{h(x)s_1(x)}\,d\mu(x),$$
which implies that 
$-\|A\|h\leq h_A\leq \|A\|h,$, for $\mu$-a.e. $x\in X$.
\par
To prove that the correspondences are inverses to each other we compute:
$$
\int_X \langle  s_1,  T P_0 A P_0 T^* s_2 \rangle d\mu = 
\langle s_1 \circ \theta_0 , P_0 A P_0 s_2 \circ \theta_0 \rangle  = \ip{s_1\circ\theta_0}{As_2\circ\theta_0}=
\int_X \langle s_1, h_0 s_2 \rangle d \mu,
$$
for every $s_1, s_2 \in E$, i.e., $h_0 = T P_0 A P_0 T^*$ and $A=A_{h_0}$.

\end{proof}

\begin{corollary}\label{corcore} Assume in addition that $h\geq c1$ for some constant $c>0$. 
Let  $T \in B( P_0 \H, \K)$ such that $T s \circ \theta_0 = h s$ for every $s \in S$. The map
$A \mapsto T P_0 A P_0 T^*$ yields a bijection between 
the commutant $C^*(X, r, m,h)'$ and the bounded harmonic functions $\mathfrak{H}:=\{f\in\M\,|\,Rf=f\}$.
\end{corollary}
\begin{proof}
Since $h\geq c1$, the boundedness condition $|h_0|\leq \text{const}\,h$ is equivalent to the essential boundedness of $h_0$, therefore it is automatically satisfied. Theorem \ref{cone} gives us the bijection between self-adjoint elements. But every element $A$ in the commutant can be written as a linear combination of self-adjoint elements: $A=(A+A^*)/2 +\text{i} (A-A^*)/2\text{i}$. Similarly for the harmonic functions. This gives the bijection.
\end{proof}

\par
Corollary \ref{corcore} shows that the set of bounded harmonic maps has also multiplicative structure, the one induced from multiplication of operators via the given bijection.
The following theorem gives an alternative, more intrinsic description of this multiplicative structure. It is also a generalization of a result from \cite{Raugi}.

\begin{theorem}\label{thprod}
Suppose $h\geq c1$ $\mu$-a.e., for some constant $c>0$. 
If $ h_1, h_2  \in  \mathfrak{H}$ then 
$$
h_1*h_2(x) := \lim_k R^k(h_1 h\- h_2)(x)
$$
exists  $\mu$ a.e and $h_1*h_2 \in \mathfrak H$. 

The map  $A \mapsto TP_0A P_0T^*$ defines a $*$-isomorphism 
from $C^*(X,r,m,h)'$ to $\mathfrak H$, when $\mathfrak H$ is equipped with  the  product $*$ and the ordinary 
involution. 
\end{theorem}

\begin{proof}
We define $\tilde R (b) = h^{-1/2}R( h^{1/2} b h^{1/2})h^{-1/2}$. Now 
  $\tilde R$ is completely positive and unital (see \cite{DuRo1}). 
  The Kadison-Schwarz inequality \cite[II.6.8.14]{MR2188261} for completely positive and unital maps yields 
  \begin{align*}
    R (b)^*h\- R(b) =  h^{1/2}\tilde R(h^{-1/2}b h^{-1/2} )^* \tilde R(h^{-1/2}bh^{-1/2})h^{1/2}  
  \leq  h^{1/2} \tilde R(h^{-1/2}b^*h\-  b h^{-1/2}) h^{1/2} = R(b^* h\- b) 
  \end{align*}  
  
  This implies that $ R(h_i^* h\- h_i ) \geq R( h_i^*) h\- R(h_i) = h_i^* h\- h_i$, $i\in\{1,2\}$,
  so $R^k(h_i^* h\- h_i)$ is a positive and increasing sequence in $\M$. Since $\sup_k \| R^k \|< \infty$, 
  the sequence is uniformly bounded, so $R^k(h_i^* h\- h_i)$ converges pointwise $\mu$-a.e. to some element $h_i*h_i$ in $\M$. And it is easy to see that $R(h_i*h_i)=h_i*h_i$.
  \par 
  Note that $h_1 h\- h_2 = \frac1 4 \sum_{k = 0}^3 i^{-k} ( h_1^* + i^k h_2)^* h\- ( h_1^* + i^k h_2)$. It follows that 
  $\lim_k R^k(h_1 h\- h_2)(x) $  exists $\mu$-a.e.  for arbitrary $h_1, h_2 \in \mathfrak H$.

If $A_1, A_2 \in C^*(X,r,m,h)'$, $h_i:= TP_0 A_i P_0 T^*$, $h_3 := TP_0 A_1 A_2 P_0 T^*$,
    then, since 
    $$\ip{s_1\circ\theta_0}{A_is_2\circ\theta_0}=\int_X\ip{s_1(x)}{h_i(x)s_2(x)}\,d\mu(x)=\ip{s_1\circ\theta_0}{(h^{-1}h_i)\circ\theta_0s_2\circ\theta_0},$$
    we see that $P_0A_iP_0$ is an operator of multiplication by $h_ih^{-1}\circ\theta_0$. 
    \par
   Let $a := h_1 h\- h_2 - h_3   $. Then
  \begin{align*}
  & \|(  P_k A_1 P_k A_2 P_k - P_k A_1 A_2 P_k) f \circ \theta_0 \|^2 =  
  \|( P_0 A_1 P_0 A_2  P_0 - P_0 A_1 A_2 P_0)U^k f \circ \theta_0\|^2 \\
  =  & \int_X \langle h\- ( h_1 h\- h_2 - h_3)m^{(k)} f \circ r^k ,
  h h\- ( h_1 h\- h_2 - h_3)m^{(k)} f \circ r^k \rangle d \mu \\
  = & \int_X \langle f , R^k (a^* h\- a) f \rangle d \mu,
  \end{align*}
  since $U$ is unitary, $P_k = U^{-k} P_0 U^k$, and $\mu$ is strongly invariant.
  \par
  Since $P_k $ converges to the identity strongly as $k$ increases, we obtain that 
  $\lim_k \int_X \langle f , R^k (a^* h\- a) f \rangle d \mu = 0$ for every $f \in \K$. But this implies that 
  there exists a subsequence such that $\lim_j R^{k_j} (a^* h\- a)(x) = 0$ $\mu$-a.e.

  The  Kadison-Schwarz inequality implies that 
  $ R^k (a)^*h\- R^k(a) \leq  R^k(a^* h\- a)$. 
  Then $R^{k_j}( a)$ converges $\mu$-a.e. to $0$, so 
  $R^{k_j}(h_1 h\- h_2)$ converges $\mu$-a.e. to $h_3$. We already know that
  $R^k(h_1 h\- h_2)$ converges $\mu$-a.e. to $h_1*h_2$, so $h_1*h_2(x) =  h_3(x)$, $\mu$-a.e.

\end{proof}
\begin{remark}\label{remprodcont}
We proved in \cite{DuRo1}, that in the case when $r$ is obtained by applying a certain covering projection to an expansive automorphism, and $m$ and $h$ are Lipschitz, with $h\geq c1$, then there is a $C^*$-algebra structure on the continuous harmonic functions. The multiplication is constructed as follows: due to the quasi-compactness of the transfer operator $R$ (restricted to Lipschitz functions), the uniform limit
$$T_1(f):=\lim_{n\rightarrow\infty}\frac{1}{n}\sum_{j=0}^{n-1}R^jf,$$
exists for every continuous $f$ and defines a continuous harmonic function. The product of two continuous harmonic functions is defined by $(h_1,h_2)\mapsto T_1(h_1h^{-1}h_2)$. Theorem \ref{thprod} shows then that this product coincides with $h_1*h_2$. In particular, if $h_1,h_2$ are continuous harmonic functions, and 
$A_{h_1},A_{h_2}$ are the associated operators in the commutant, then $A_{h_1}A_{h_2}=A_{h_1*h_2}$ is also associated to a continuous harmonic function $h_1*h_2=T_1(h_1h^{-1}h_2)$.
\end{remark}

The next corollary shows how the covariant representation can be decomposed using projections in the algebra $\mathfrak H$.

\begin{corollary}
  Suppose $h_1, \dots, h_l$ form a family of mutually  orthogonal projections in $\mathfrak H$, i.e., 
  \begin{align*}
    h_i * h_i & = h_i^* = h_i ,~ 1 \leq i \leq l,\quad
    h_i * h_j  = \delta_{i,j}, ~1 \leq i,j \leq l.
  \end{align*}
  Let $p_1, \dots, p_l$
  denote the corresponding orthogonal projections  in  $C^*(X, r, m,h)'$.
   Moreover, let $(  \H_i, \langle \cdot , \cdot \rangle_i)$ denote 
  the Hilbert space,  associated to
  $h_i$,   in  Theorem \ref{representation} with  the representation 
  $\pi_i : C(X) \rightarrow B(\H_i)$, unitary $U_i$ and 
  increasing family of projections $\{P_k^i\}_k$.

  Then $p_i \H = \H_i$, $\langle  p_is_1, p_is_2 \rangle =  
  \langle s_1, s_2 \rangle_i$
  for every $s_1, s_2$ sections in $\theta_0^*\xi$.
  We obtain an isometry $$J : \oplus_{i=1}^l \H_i \rightarrow \H,\quad J(s_1,\dots,s_l)=p_1s_1+\dots p_ls_l,$$
  such that 
  \begin{enumerate}
  \item $\pi(a) J = J ( \oplus_{i=1}^l \pi_i(a))$
  \item $UJ = J ( \oplus_{i=1}^l U_i)$
  \end{enumerate}
  If $\sum_{i=1}^l h_i = h$ then $J$ is a unitary.

\end{corollary}

\section{Cocycles}\label{cocy}
In this section we give an alternative description of the Hilbert space $\H$ and of the operators in the commutant in terms of some matrix valued measures $P_x$ on the solenoid $X_\infty$.

In the scalar case, the measures $P_x$ are random walk measures with variable coefficients. Each point in the solenoid $X_\infty$ can be constructed in the following way. Pick a point $x_0\in X$. Then, since $rx_1=x_0$, one has to make a choice $x_1$ out of the finitely many roots in $r^{-1}x_0$. The QMF equation amounts to
$$\frac{1}{\# r^{-1}x}\sum_{ry=x}|m_0(y)|^2=1.$$
Therefore $|m_0(y)|^2/\#r^{-1}(y)$ can be interpreted as the probability of transition from $x$ to its root $y$. 
At the next step one makes a transition from $x_1$ to its root $x_2$ with probability given by $|m_0(x_2)|^2/\#r^{-1}x_1$. And so on. The measure $P_x$ is the path measure obtained in this way.

In our matricial case, the measures $P_x$ will be operator valued, but the idea stays the same. However some complications appear because we are dealing with a non-commutative situation.
 
 We will also see that the operators in the commutant of the covariant representation are in fact multiplication operators by matrix valued functions. This will enable us to give a more concrete form of the correspondence between harmonic maps and the operators in the commutant. The result resembles the Poisson-Fatou-Privalov theorems in harmonic analysis: the harmonic map is the integral on the boundary of the operator in the commutant, and the operator is a radial limit of the harmonic map. In our case the boundary is the solenoid $X_\infty$.

First, by applying the Trace we will convert the matrix measures into a scalar measure and a matrix valued Radon-Nykodim derivative.
\par
Let $ C_k \subset C(X_\infty) $ be  the set of continuous functions on $X_\infty$ that only depend on the $k+1$ first 
coordinates. Then
$$f \circ \theta_k \mapsto \int_X  \text{Tr}( R^k ( h f)(x) )   d \mu(x) $$
defines a positive linear functional $\tau_k \in C_k^*$ such that the compatibility condition
$\tau_{k+1} | _{C_k} = \tau_k$ holds.
\par 
Since $\sup_k \| R^k\| < \infty$  and $\cup_k C_k$ is dense in $C(X_\infty)$,  there 
exists a positive functional $\tau \in C(X_\infty)^*$ such that 
$\tau|_{C_k} =\tau_k$ for every $k$. 
Let $\hat \mu$ be the 
measure on $X_\infty$, provided by the Riezs representation theorem 
such that $\tau(f) = \int_{X_\infty} f d \hat \mu$

\begin{proposition}\label{prop4-3}{\rm (i)}
  There exist a positive $\Delta(x) \in \text{End}_\C( \theta_0^* \xi |_x )$, for $\hat \mu$-a.e. $x\in X_\infty$ such that $\| \Delta(\cdot)\| \in L^\infty(X_\infty, \hat \mu)$
  and 
  $$
  \int_{X_\infty} \langle s_1 (x), \Delta(x) s_2(x) \rangle d\hat \mu = \langle s_1, s_2 \rangle,
  $$
  for every pair of sections $s_1, s_2 : X_\infty \rightarrow \theta_0^* \xi$.
  \par
  {\rm (ii)} If $A$ is an operator in the commutant $C^*(X,r,m,h)'$ then for $\hat\mu$-a.e. $x\in X_\infty$, there exists $A(x)\in\mbox{End}_\C(\theta_0^*\xi|_x)$, such that 
  $As(x) = A(x) s(x)$, $\hat\mu$-a.e., and 
  \begin{equation}\label{eqcocy} \Delta(x) m \circ \theta_0 (x) A (\hat r (x))  = \Delta(x) A(x) m \circ \theta_0(x),\hat\mu\mbox{-a.e.}\end{equation}
  Also, conversely, any such esentially bounded matrix-valued function $x\mapsto A(x)$, satisfying \eqref{eqcocy}, defines an operator in the commutant $C^*(X,r,m,h)'$. 
  \end{proposition}

  \begin{proof}
    Denote  $L^2( \theta_0^* \xi, \hat \mu)$ denote the $L^2$-sections in the bundle $\theta_0^*\xi$ with respect to the measure $\hat \mu$.
    Note that if $s$ is a section in $\xi$, then 
    \begin{align*}
      \langle s \circ\theta_k , s \circ \theta_k \rangle & = \int_X \sum_{  r^k y = x} \frac{1}{\#r^{-k}(x)}
      \langle m^{(k)}(y) s(y)  , h(y) m^{(k)}(y) s(y) \rangle d \mu \\
      & \leq \int_X \sum_{  r^k y = x} \frac{1}{\#r^{-k}(x)} \|  m{(k)^*}(y)h(y) m^{(k)}(y)\| \langle s(y), s(y) \rangle d \mu(x) \\
      & \leq \int_X  \sum_{ r^k y = x} \frac{1}{\#r^{-k}(x)} 
      \text{Tr}  (m^{(k)*}(y) h(y) m^{(k)}(y))   \langle s(y), s(y) \rangle d \mu(x) \\
      & = \tau_k ( \langle s \circ \theta_k , s \circ \theta_k \rangle ) \\
      & = \int_{X_\infty} \langle  s \circ \theta_k , s \circ \theta_k \rangle d\hat \mu,
    \end{align*}
    so there exists a unique and positive bounded operator on $L^2( \theta_0^* \xi, \hat \mu)$ 
    such that 
     $
  \int_{X_\infty} \langle s_1 (x), (\Delta s_2)(x) \rangle d\hat \mu = \langle s_1, s_2 \rangle, 
  $
    for every pair of sections $s_1, s_2 : X_\infty \rightarrow \theta_0^* \xi$.

    The $C(X_\infty)$ module structure on the section in $\theta_0^* \xi$ gives us representations 
    of $C(X_\infty)$ on $\H$ and $L^2 ( \theta_0^* \xi, \hat \mu)$ as follows: 
    Every $f \in C(X_\infty)$ gives a multiplication operator $M_f \in B(\H)$
    such that $(M_f s)(x) = f(x) s(x)$. We see that $M_f^* = M_{\overline f}$
    as operators on both $\H$ and $L^2(\theta_0^* \xi, \hat \mu)$, so 
    $$
      \int_{X_\infty} \langle s_1 (x), (\Delta M_f s_2)(x) \rangle d\hat \mu  = 
     \langle s_1, M_f s_2 \rangle  =  \langle M_{\overline f} s_1,  s_2 \rangle $$
       $$=  \int_{X_\infty} \langle M_{\overline f} s_1 (x), \Delta  s_2(x) \rangle d\hat \mu
      =  \int_{X_\infty} \langle s_1 (x), M_f \Delta  s_2(x) \rangle d\hat \mu,$$
     and therefore $M_f \Delta = \Delta M_f$. This implies that $\Delta s(x) = \Delta(x) s(x)$
     for a $\Delta(x) \in \text{End}(\theta_0^* \xi|_x)$ as described, for $\hat \mu$-a.e. $x\in X$.
  \par
  (ii) First note that $C^*( X, r, R, h)' \subset C(X_\infty)'$. This is because $U^{-1}\pi(f)U=M_{f\circ\theta_1}$, so if an operator commutes with $U$ and $\pi$ then it must commute with all multiplications by functions which depend only on finitely many coordinates. But these functions are dense in $C(X_\infty)$ and we obtain the inclusion.
  \par
  From this, we see that there exist $A(x) \in \theta_0^* \xi |_x$ such that 
$As(x) = A(x) s(x)$, $\hat \mu$-a.e. Moreover, $U A = A U$ implies the relation 
\eqref{eqcocy}. The converse follows by a computation to prove that $A$ defined from $x\mapsto A(x)$ commutes with $U$ and $\pi(f)$ for all $f\in C(X)$.
\end{proof}

\begin{definition}
Because of the relation \eqref{eqcocy}, if $A$ is an operator in the commutant $C^*(X,r,m,h)'$ then we call $A$ a {\it cocycle}. 
\end{definition}
\par
In the next proposition we will define the positive-matrix-valued measures $P_x$ which can be used to represent the inner-product on the Hilbert space $\H$. 
\begin{proposition}\label{proppx}
For each $x\in X$, let $\Omega_x = \{ y \in X_\infty | \theta_0 (y) = x\}$. 
  There exists  a positive operator-valued  measure $P_x$ from the Borel sets in  $\Omega_x$ 
  to  $\text{End}_\C( \xi|_x)$ such that 
  % $x \mapsto \int_{\Omega_x} f d P_x$ is measurable 
  $\int_{\Omega_x} f \circ \theta_k dP_x = R^k ( f h)(x),$
  for every bounded measurable function $f$ on $X$,
  and 
  $$
  \int_X \int_{\Omega_x} \langle s_1 (y),  dP_x(y) s_2(y) \rangle d \mu(x) = 
  \langle s_1 , s_2 \rangle,
  $$
  for every pair of sections $s_1, s_2 : X_\infty \rightarrow  \theta_0^* \xi$.
\end{proposition}
\begin{proof}
  Note that $s_i(y) \in \theta_0^* \xi|_x$ for every $y \in \Omega_x$. 
  Let $C_{k,x} \subset C(\Omega_x)$ denote the functions that only depend on the first $k+1$
  variables. Moreover, define 
  $\sigma^x_k : C(\Omega_x) \rightarrow \text{End}( \xi|_x)$ by
  $$\sigma_k^x ( f \circ \theta_k) = R^k( f h)(x)$$
  We see that 
  $\sigma_k^x$ defines a positive bounded operator from $C_{k,x}$ to $M_d(\C)$. Moreover the compatibility condition 
  $\sigma_{k+1}^x |_{ C_k}= \sigma_k^x$ holds. And, since    $\sup_k  \| R^k\| < \infty$ and 
  $\cup_k C_{k,x} \subset C(\Omega_x)$ is dense, we see that there exists a 
  positive and bounded linear map $P_x : C(\Omega_x) \rightarrow \text{End}(\xi|_x)$
  for almost every $x \in X$, such that $\sigma^x |_{C_{k,x}} = \sigma_k^x$. A matrix computation 
  implies that 
   $$
  \int_X \int_{\Omega_x} \langle s_1\circ \theta_k (y),  dP_x(y) s_2\circ \theta_k (y) \rangle d \mu(x) = 
  \langle s_1 \circ \theta_k  , s_2 \circ \theta_k \rangle,
  $$
  for every pair of sections $s_1, s_2 : X \rightarrow \xi$. The last claim follows from this by a density argument. 
\end{proof}
The next Lemma can be obtained by direct computation.
\begin{lemma}\label{measures}
For $x\in X$, let $\hat\mu_x$ be the measure on $\Omega_x$ defined by: $$\int_{\Omega_x}f\circ\theta_k\,d\hat\mu_x=\Tr(R^k(fh)(x))=\Tr(\int_{\Omega_x}f\circ\theta_k\,dP_x).$$
Then, for all bounded measurable functions on $X_\infty$
$$\int_{X_\infty}f\,d\hat\mu=\int_X\int_{\Omega_x}f\,d\hat\mu_x\,d\mu(x).$$
If $C_{\omega_1,\dots,\omega_n}$ is the set of points $(x,x_1,x_2,\dots)\in\Omega_x$ such that $x_1=\omega_1,\dots,x_n=\omega_n$, then
$$P_x(C_{\omega_1,\dots,\omega_n})=\frac{1}{\#r^{-n}(x)}({m^{(n)}}^*hm^{(n)})(\omega_n),\quad
\hat\mu_x(C_{\omega_1,\dots,\omega_n})=\frac{1}{\#r^{-n}(x)}\Tr\left(({m^{(n)}}^*hm^{(n)})(\omega_n)\right).$$
\end{lemma}

Having defined the matrix valued measures $P_x$, the correspondence between cocyles and harmonic functions in Theorem \ref{cone}(ii) can be given now in terms of a matrix valued conditional expectation:
\begin{proposition}\label{expe}
  Let $A \in C^*(X, r, m, h)'$ and $A(x) \in \text{End}( \theta_0^* \xi |_x)$ such that 
  $As(x) = A(x) s(x)$, $\hat \mu$-a.e., for every section $s$. We have the following identity:
  $$
  TP_0 A P_0 T^* (x) = \int_{\Omega_x} dP_x(y)  A(y),
  $$
  i.e., $x \mapsto \int_{\Omega_x} dP_x(y)  A(y) $ is a fixed point for  $R$.
  \begin{proof}
    If  $s_1, s_2:  X \rightarrow \xi$ are sections, then  
  \begin{align*}
     & \int_X \langle s_1(x) , (T P_0 A P_0 T^*)(x) s_2(x) \rangle d \mu(x)  = 
    \langle s_1 \circ \theta_0 , P_0 A P_0 s_2 \circ \theta_0  \rangle \\ = &  
    \langle s_1\circ \theta_0, A s_2 \circ \theta_0 \rangle 
     = \int_X \langle s_1(x) , \int_{\Omega_x} d P_x(y) A(y) s_2(x) \rangle d \mu(x).
  \end{align*}
  \end{proof}
\end{proposition}
For the inverse correspondence in Theorem \ref{cone}(i), from harmonic functions to cocycles, we have the following result:

\begin{theorem}\label{conv} Assume in addition that $h\geq c1$,
$\mu$-a.e., for some constant $c>0$.
Let $h_0$ be a bounded harmonic function and let $A$ be the corresponding
cocycle as in Theorem \ref{cone}(i) and Proposition \ref{prop4-3}. Then
$$\lim_{k\rightarrow\infty}\frac{({m^{(k)}}^*h_0m^{(k)})\circ\theta_k}{\Tr({({m^{(k)}}^*hm^{(k)})\circ\theta_k})}=\Delta
A,$$
pointwise $\hat\mu$-a.e.
\end{theorem}
\begin{proof}
For $f$ bounded measurable function on $X_\infty$ let $E_k(f)$ denote the conditional expectation onto the functions that depend only on the first $k+1$ coordinates, with respect to the measure $\hat\mu$. If $F$ is a matrix-valued function on $X_\infty$, then $E_k(F)$ is the matrix-valued function obtained by applying $E_k$ to each component.
\par
Let $A_k:=({m^{(k)}}^{-1}h^{-1}h_0m^{(k)})\circ\theta_k$. (Recall that $m(x)$ is invertible for $\mu$-a.e. $x\in X$). Then 
$$\ip{f\circ\theta_k}{A_kg\circ\theta_k}=\ip{f\circ\theta_k}{Ag\circ\theta_k}=\ip{f\circ\theta_k}{P_kAP_kg\circ\theta_k}.$$

Take $f_k=f\circ\theta_k,g_k=g\circ\theta_k$. 
$$\int_X\int_{\Omega_x}\ip{f_k}{E_{k+1}(\Delta)A_{k+1}g_k}\,d\hat\mu_x\,d\mu(x)=
\int_X\int_{\Omega_x}\ip{f_k}{\Delta A_{k+1}g_k}\,d\hat\mu_x\,d\mu(x)$$
$$
=\ip{f_k}{A_{k+1}g}=\ip{f_k}{P_{k+1}AP_{k+1}g_k}=\ip{f_k}{Ag_k}=\ip{f_k}{A_kg_k}$$
$$=\int_X\int_{\Omega_x}\ip{f_k}{\Delta A_kg_k}\,d\hat\mu_x\,d\mu(x)
=\int_X\int_{\Omega_x}\ip{f_k}{E_k(\Delta)A_kg_k}\,d\hat\mu_x\,d\mu(x).$$
So for $\mu$-a.e. $x\in X$, $E_k(E_{k+1}(\Delta)A_{k+1})(x,\cdot)=(E_k(\Delta)A_k)(x,\cdot)$, $\hat\mu_x$-a.e., for all $k\geq1$. Therefore the sequence $\{(E_k(\Delta)A_k)(x,\cdot)\}_k$ is a martingale. By Doob's martingale convergence theorem, 
$E_k(\Delta)A_k(x,\cdot)$ converges pointwise $\hat\mu_x$-a.e. to $\Delta A(x,\cdot)$.

\par 
Now we compute $E_k(\Delta)(x,\cdot)$. We have
\begin{equation}\label{eqekdelta0}
\int_{\Omega_x}\ip{f_k}{E_k(\Delta)g_k}\,d\hat\mu_x=\int_{\Omega_x}\ip{f_k}{\Delta g_k}\,d\hat\mu_x=\frac{1}{\#r^{-k}(x)}\sum_{r^ky=x}\ip{f(y)}{{m^{(k)}}^*h(y)m^{(k)}(y)g(y)}.
\end{equation}
Let $C_{\omega_1,...,\omega_k}$ be the cylinder of points in $\Omega_x$ that start with $\omega_1,...,\omega_k$. Let $f_k$ and $g_k$ be supported on $C_{\omega_1,...,\omega_k}$.
Then we obtain from \eqref{eqekdelta0} that
$$E_k(\Delta)(x,\omega_1,\dots,\omega_k)\hat\mu_x(C_{\omega_1,...,\omega_k})=\frac{1}{\#r^{-k}(x)}\left({m^{(k)}}^*hm^{k}\right)(\omega_k).$$
But,  with Lemma \ref{measures}, 
$$\hat\mu_x(C_{\omega_1,...,\omega_k})=\Tr\left(P_x(C_{\omega_1,...,\omega_k})\right)=\operatorname*{Tr}\left(\frac{1}{\#r^{-k}(x)}({m^{(k)}}^*hm^{(k)})(\omega_k)\right).$$
Thus
\begin{equation}\label{eqe_kdelta}
E_k(\Delta)=\frac{({m^{(k)}}^*hm^{(k)})\circ\theta_k}{\Tr(({m^{(k)}}^*hm^{(k)})\circ\theta_k)}
\end{equation}
Then
$$E_k(\Delta)A_k=\frac{({m^{(k)}}^*h_0m^{(k)})\circ\theta_k}{\Tr({({m^{(k)}}^*hm^{(k)})\circ\theta_k})}.$$
This proves the theorem.
\end{proof}

\begin{remark}
In the scalar case, the terms involving $m^{(n)}$ will disappear and we reobtain the results from \cite{DuJomart}. 
\end{remark}

\begin{remark}{\bf An ergodic limit for low-pass filters.} Our limit theorem can be used to obtain an interesting limit result. The matrices enable us to use a trick to compare two different measures. 

Let $X=\bt^1$, $r(z)=z^2$. Denote by $\psi_1,\psi_2$ the inverse branches of $r$, $\psi_1(e^{i\theta})=e^{i\theta/2}$, $\psi_1(e^{i\theta})=e^{i(\theta+2\pi)/2}$, for $\theta\in[0,2\pi)$. Then the solenoid $X_\infty$ is measurably isomorphic to $X\times\Omega$, where $\Omega:=\{1,2\}^{\bn}$, because a point $(z_0,z_1,\dots)\in X_\infty$ consists of $z_0\in\bt^1$ and a choice of the inverse branches $\omega_1,\omega_2,\dots$. \par
Let $m_1(z)=(1+z)/\sqrt{2}$ (or any low-pass filter that gives orthogonal scaling functions in $L^2(\br)$), and $m_2(z)=1$. Let $h_1=h_2=1$. Then the covariant representation for $(m_1,h_1)$ is $L^2(\br)$ so the measure $P_x^1$ on the solenoid is supported on sequences $\omega\in\Omega$ such that there exists $n_0$ such that $\omega_n=0$ for all $n\geq n_0$ or $\omega_n=1$ for all $n\geq n_0$ (see \cite{DuJoifs}, or Section \ref{low}). The covariant representation for $(m_2,h_2)$ is on the solenoid $X_\infty$ with the Haar measure $\mu_H$, and the measures $P_x^2$ are the Bernoulli measures on $\Omega$ where $\{0\}$ and $\{1\}$ have equal probabilities $1/2$ (see also \cite{DuJoifs} and \cite{DuJowf}). 
\par
Let $m:=\left[\begin{array}{cc}m_1&0\\ 0&m_2\end{array}\right].$ Then the measure $P_x$ is clearly 
$\left[\begin{array}{cc}P_x^1&0\\ 0&P_x^2\end{array}\right]$, i.e., for $f\in C(X_\infty)$,
$$\int_{X_\infty}f\,d\,P_x=\left[\begin{array}{cc}\int_{X_\infty}f\,dP_x^1\,d\mu(x)&0\\ 0&\int_{X_\infty}f\,dP_x^2\,d\mu(x)\end{array}\right].$$
The trace of this measure is $\hat\mu_x=P_x^1+P_x^2$. As we explained before (see also Section \ref{low}), the measure $P_x^1$ is atomic, and let us denote the support of this measure by $R_x$. On the other hand, $P_x^2$ is the Bernoulli measure, so it has no atoms, therefore $R_x$ has $P_x^2$-measure zero. \par
Then it is easy to see that  $\frac{dP_x^1}{d\hat\mu_x}=\chi_{R_x}$ and $\frac{dP_x^2}{d\hat\mu_x}=\chi_{\Omega\setminus R_x}$. Then 
$$\Delta=\frac{dP_x}{d\hat\mu_x}=\left[\begin{array}{cc}\chi_{R_x}&0\\ 0&\chi_{\Omega\setminus R_x}\end{array}\right].$$
Note that this shows that $\Delta(x)$ is singular everywhere. 
\par
We will remark that this implies an interesting phenomenon, which occurs for any low-pass filter that gives orthogonal scaling functions in $L^2(\br)$. 
\par
First, as shown in \cite{DuJoifs}, since $P_x^1$ is supported in on $R_x$, it follows that for any $\omega$ outside $R_x$ one has $$\lim_{k\rightarrow\infty}\frac{|m_1^{(k)}(\psi_{\omega_k}\dots\psi_{\omega_1}x))|^2}{2^k}=P_x^1(\{\omega\})=0.$$
\par
However, from the convergence Theorem \ref{conv}, we have that for $\mu$-a.e. $x\in\bt^1$ and for $\hat\mu_x$-a.e. $\omega$:
$$\lim_{k\rightarrow\infty}\frac{{m^{(k)}}^*m^{(k)}(\psi_{\omega_k}...\psi_{\omega_1}x)}{\Tr({m^{(k)}}^*m^{(k)}(\psi_{\omega_k}...\psi_{\omega_1}x))}=\Delta(x,\omega).$$
This implies that for $\hat\mu_x$ a.e. $\omega\in \Omega\setminus R_x$
$$0=\lim_{k\rightarrow\infty}\frac{|m_1^{(k)}(\psi_{\omega_k}...\psi_{\omega_1}x)|^2}{|m_1^{(k)}(\psi_{\omega_k}...\psi_{\omega_1}x)|^2+1},$$
so $\lim_{k\rightarrow\infty}m^{(k)}(\psi_{\omega_k}...\psi_{\omega_1}x)=0$ for $\hat\mu_x$-a.e. $\omega$, so
we obtain the much stronger limit for $P_x^2$-a.e. $\omega$ (and recall that $P_x^2$ is the Bernoulli measure):
$$\lim_{k\rightarrow\infty}m^{(k)}(\psi_{\omega_k}...\psi_{\omega_1}x)=0.$$
\end{remark}

\section{Low-pass filters}\label{low}
In the scalar case, if the filter $m_0$ satisfies a low-pass condition $m_0(1)=\sqrt{2}$, the classical wavelet theory shows that the scaling equation has a solution in $L^2(\br)$. Of course the solution might be a non-orthogonal scaling function, and then there are super-wavelet constructions (see \cite{BDP}) that will give orthogonal solutions. functions. 

The point we want to make is that, when a low-pass condition is satisfied, the resulting covariant representation can be realized on $L^2(\br)$, as in classical wavelet theory, or in a direct sum of copies of $L^2(\br)$, as in the super-wavelet theory developed in \cite{BDP}. 

In the matricial case, the low-pass condition is replaced by the $E(l)$-condition introduced in \cite{Jiang}.

Our covariant representations are on the solenoid $X_\infty$, but we show that when a low-pass condition is satisfied the measures $P_x$ are supported on an embedding of $\br$ in the solenoid, and $P_x$ is directly related to the scaling functions.

We recall now the setup from \cite{DuRo1}. 
We assume that $\tilde X$ is a complete metric space with an isometric  covering space 
group action of a group $G$ such that $\tilde X /G =: X$ is compact. Moreover, we assume that $\tilde r : \tilde X \rightarrow \tilde X$ is a strictly expansive homeomorphism and there exists an endomorphism $A\in \text{End}(G)$ such that 
$ \tilde r g = (A g)\tilde r  $ and $AG$ is a normal subgroup of index $q$. Since $\tilde r$ is expansive it has a fixed point $\tilde x_0$.
\par
Let $p : \tilde X \rightarrow X$ denote the quotient covering map and define $r : X \rightarrow X$ as
$r(p(x)) = p( \tilde rx)$. Let $x_0:=p(\tilde x_0)$. Let $\mu$ be a strongly invariant measure on $X$ and $\tilde\mu$ the measure on $\tilde X$ obtained by lifting the measure $\mu$ by the covering map $p$ (see also \eqref{mutilde}).
\par
We assume that the bundle $\xi$ over $X$ is a Lipschitz continuous bundle, and 
that $p^* \xi$ is trivial.
\begin{example}\label{mainex}
The main example is the one used in wavelet theory: $\tilde X=\br^n$. The group $G=\bz^n$ acts on $\br^n$ by translations: $g,x\mapsto x+g$, ($x\in\br^n, g\in\bz^n$). Define the map $\tilde r(x)=Ax$, where $A$ is an $n\times n$ expansive integer matrix. The quotient $\br^n/\bz^n$ can be identified with the torus $\bt^n$, $p:\br^n\rightarrow\br^n/\bz^n$ is the quotient map, and let $r(x)=Ax\mod\bz^n$ for $x\in\bt^n$. The fixed point of $\tilde r$ is $\tilde x_0=0$, and $x_0:=p(\tilde x_0)=0$. The bundle $\xi$ over $\bt^n$ is $\bt^n\times\bc^d$, and $p^*\xi=\br^n\times\bc^d$.
\end{example}
\par 
We let $m \in M_d(\text{Lip}_1(X))$.
Moreover, we assume that 
\begin{itemize}
\item $R1 = 1$;
\item The matrix $m(x_0)/\sqrt{q}$ satisfies the $E(l)$ condition (according to \cite{Jiang}), i.e., $1$ is the only eigenvalue of $m(x_0)/\sqrt{q}$ with absolute value greater than or equal to 
$1$, and  its algebraic and geometric multiplicity are both equal to $l\geq 1$. Let $E_1$ denote the eigenspace corresponding to the eigenvalue $1$.
\end{itemize}
\par
Recall that $S$ is the set of continuous sections in the bundle $\xi$.\par Let $\Xi := \{  f \in C_b(\tilde X) | \sum_{g \in G} |f |^2 \circ g \in C(X) \} $. 
 $\Xi$ is a $C(X)$-Hilbert module with the inner-product: 
 $$
 \langle \zeta, \eta \rangle' = \sum_{g \in G} (\overline \zeta \eta ) \circ g.
 $$
This module inner-product is related to the ordinary inner-product from $L^2(\tilde X, \tilde \mu)$, by
\begin{equation}\label{mutilde}
\int_{\tilde X} \overline \zeta \eta d \tilde \mu = \int_X \langle \zeta, \eta \rangle' d \mu,\quad(\zeta,\eta\in\Xi). 
\end{equation}
 Let $\tilde  U \in B(L^2(\tilde X , \tilde \mu))$ denote  the unitary operator 
 defined by $$\tilde U f = q^{1/2} f \circ \tilde r,\mbox{ and }, \tilde U_l:=\underbrace{\tilde U\oplus\dots\oplus\tilde U}_{l\mbox{ times }}.$$ Define the representation $\tilde \pi_l : C(X) \rightarrow B( \oplus_{j=1}^lL^2(\tilde X, \tilde \mu))$
by $$\tilde \pi_l(a) (f_1, \dots, f_l)(x) = ( a(px)f_1(x), \dots, a(px)f_l(x)),\quad(a\in C(X),f_1,\dots,f_l\in L^2(\tilde X,\tilde\mu)).$$
\par
We will need to define some ``scaling functions''. These will be fixed points of a refinement operator obtained as the limits of the iterates of this refinement operator. \par
To initialize the iteration, we fix an $f \in \Xi$ such that:
\begin{itemize}
\item $\langle f, f \rangle' = 1$
\item $( \sum_{g \in G} | f(gx) -  f(gy)|^2)^{1/2} \leq D d(x,y)$
\item $f(gx_0) = 0 $ for every $g \neq 1$.
\end{itemize}
\par
Let 
$s_1, \dots, s_l\in S$ be Lipschitz sections in $\xi$ such that $s_1(x_0), \dots , s_l(x_0)$ form an orthonormal basis for 
 $ E_1 $. 
\par
As in \cite{DuRo1}, we define the starting points for the cascade algorithm $W_j \in \text{Hom}_{C(X)}(S, \Xi)$ by 
$W_j s = \langle s_j , s \rangle f$ and the {\it refinement operator} $$M : \text{Hom}_{C(X)}(S, \Xi) \rightarrow \text{Hom}_{C(X)}(S, \Xi),\quad
(M W) s = \tilde U ^{-1} Wm s \circ r ,\mbox{ for }W \in  \text{Hom}_{C(X)}(S, \Xi) ,s \in S.$$
\par

\begin{proposition}\label{duro}\cite{DuRo1}
\begin{enumerate}
\item The following limit exists and is uniform on compact sets:
$$\mathcal P(\tilde x):=\lim_{k\rightarrow\infty}q^{-k/2}m^{(k)}(p(\tilde r^{-k}\tilde x)),\quad(\tilde x\in \tilde X).$$
Also, $\mathcal P(\tilde x_0)$ is the projection onto the eigenspace $E_1$.
\item For each $j\in\{1,\dots,l\}$ and for $s\in S$, $\{M^kW_js\}_{k\geq 1}$ converges uniformly on compact sets to $\mathcal W_js$, and this defines $\W_j\in \text{Hom}_{C(X)}(S, \Xi)$, $\W_js(\tilde x)=\ip{s_j(x_0)}{\mathcal P(\tilde x)s(p\tilde x)}$, $s\in S$, $\tilde x\in\tilde X$. Moreover $M \W_j = \W_j$. 
\item For each $j\in\{1,\dots,l\}$ the map $h_j:=\W_j^*\W_j$ defines a minimal projection in the algebra of continuous harmonic functions $\mathfrak H_c$ ($h_j$ is a projection also in the algebra of bounded measurable harmonic functions $\mathfrak H$ but it is not necessarily minimal).
$$h_j(px)s(px)=\sum_{g\in G}\ip{\mathcal P^*(gx)s_j(x_0)}{s(px)}\mathcal P^*(gx)s_j(x_0),\quad (s\in S,x\in\tilde X).$$
Moreover these projections are mutually orthogonal in this algebra. 
\end{enumerate}

\end{proposition}

\begin{theorem}\label{thj}
Let $h := \sum_{j=1}^l h_j$ and let $\H$ denote Hilbert space of  
the covariant representation  obtained from $m$ and $h$. There exists a unitary $J : \H \rightarrow \bigoplus_{j=1}^l L^2(\tilde X, \tilde \mu)$ 
such that 
\begin{enumerate}
\item 
  $\tilde U_l J = JU$
\item 
  $J s \circ \theta_0 = \bigoplus_{j=1}^l \W_j s$ for every $s \in S$.
\item 
  $\tilde \pi_l(a) J = J \pi(a)$ for every $a \in C(X)$.
\end{enumerate}
\end{theorem}

\begin{proof}
 Define $J_k  : \H_k  \rightarrow \bigoplus_j L^2(\tilde X , \mu)$ by
$J_k U^{-k} s \circ \theta_0  = \bigoplus_j \tilde U^{-k} \W_j s   $.
$J_k$ is an isometry because
\begin{align*}
   \langle J_k &U^{-k} s \circ \theta_0 , J_k U^{-k} s \circ \theta_0 \rangle_{L^2} = 
  \sum_j  \langle \tilde U ^{-k} \W_j s ,  \tilde U ^{-k} \W_j s \rangle_{L^2} 
   = \sum_j \langle \W_j s, \W_j s \rangle_{L^2} 
   =    \sum_j \int_X \langle \W_j s, \W_j s \rangle' d \mu \\
   = &  \int_X \langle s, \sum_j \W_j^*\W_j s \rangle d\mu =\int_X\ip{s}{hs}\,d\mu 
   = \langle s \circ \theta_0 , s \circ \theta_0 \rangle 
   =   \langle U^{-k}s \circ \theta_0 , U^{-k} s \circ \theta_0 \rangle 
\end{align*}

Moreover, since $\W_j=M\W_j=\tilde U^{-1}\W_jms\circ r$,
$$J_k  U^{-(k-1)}  s \circ \theta_0 = J_kU^{-k}Us\circ\theta_0=\oplus_j \tilde U^{-k}\W_j m s\circ r = \oplus_j \tilde U^{-(k-1)} \W_j s
= J_{k-1}U^{-(k-1)}s \circ \theta_0$$
so $J_k |_{\H_{k-1}} = J_{k-1}$.
This gives us an isometric map $J : \H \rightarrow \bigoplus _j L^2( \tilde X, \tilde\mu)$, and the intertwining properties of $J$ are checked by a direct computation.

Let $\tilde \H := \cup_k J \H_k \subset \bigoplus_j L^2(\tilde X, \tilde \mu)$ and let $Q$ denote the orthogonal projection onto this space. 
We have $\tilde U_l Q = Q \tilde U_l$,  
$\tilde \pi_l (a) Q = Q \tilde \pi_l(a)$
for every $a \in C(X)$.
\par
We have that $\tilde U_l^{-k}\tilde \pi_l(a)\tilde U_l^k$ is a multiplication by the matrix that has $a\circ p\circ\tilde r^{-k}$ on the diagonal. Since $\tilde r$ is expansive, an application of the Stone-Weierstrass theorem shows that $\{a\circ p\circ \tilde r^{-k}\,|\,a\in C(X),k\geq 0\}$ is dense in $C_c(\tilde X)$. Since $Q$ commutes with all operators of the form $\tilde U_l^{-k}\tilde \pi_l(a)\tilde U_l^k$, this implies that $Q$ commutes with $L^\infty(\tilde X,\tilde \mu)\otimes I_l$. Since $L^\infty(\tilde X,\tilde\mu)$ is a maximal abelian subalgebra, it follows that $Q$ corresponds to a pointwise multiplication by a map in $L^\infty( \tilde X, \tilde \mu)\otimes M_l(\bc)  $.
Since $Q$ is a projection, this  map is projection valued. Moreover, $Q(x) = Q(\tilde rx)$ 
$\tilde \mu$-a.e., since $Q\tilde U_l = \tilde U_l Q$.
\par
Consider now $\varphi_i:=\oplus_{j=1}^n \W_js_i$ for $i\in\{1,\dots,l\}$. We have with Proposition \ref{duro}, $\varphi_i(\tilde x_0)=e_i$, the canonical vectors in $\C^l$. Then, using the continuity of $\W_j$, we have that 
for $\tilde x$ in a neighborhood of $\tilde x_0$, $\{\varphi_i(\tilde x)\,|\,i\in\{1,\dots,l\}\}$ forms a basis for $\C^l$. Since $\varphi_i\in \tilde \H$, $Q\varphi_i=\varphi_i$ so $Q(x)\varphi_i(x)=\varphi_i(x)$ for all $x\in \tilde X$. But then, $Q(x)$ must be the identity in a neighborhood of $\tilde x_0$, and since $Q(x)=Q(\tilde rx)$, and $\tilde r$ is expansive, we obtain that $Q(x)$ is the identity for all $x\in \tilde X$. Thus $Q=1$ and $\tilde\H$ is the entire space $\bigoplus_j L^2(\tilde X,\tilde \mu)$.

\end{proof}

\begin{proposition}
The map $\hat i:\tilde X\rightarrow X_\infty$, $\hat i(x)=(p(\tilde r^{-k}x))_{k=1}^\infty$ is a continuous bijection onto the set of sequences $(z_k)_k\in X_\infty$ with $\lim_{k\rightarrow\infty}z_k=x_0$. For all $p(x)\in X$, $\hat i(\tilde X)\cap \Omega_{px}=\hat i(Gx)$.
\par
Let $P_x$ be the measures associated to $h=\sum_j h_j$ as in Proposition \ref{proppx}. Then $P_{px}$ is atomic and supported on $\hat i(Gx)$, and 
$$P_{px}(\{\hat i(gx)\})=\mathcal P(gx)^*\mathcal P(gx),\quad(x\in\tilde X).$$

\end{proposition}
\begin{proof}
Since $\tilde r$ is expansive, the sequence $\tilde r^{-k}x$ converges to the fixed point $\tilde x_0$, so for $k$ large $\tilde r^{-k}x$ is in some neighborhood where $p$ is injective. This implies that $\hat i$ is injective. 
\par
The continuity of $\hat i$ is clear, and $p(\tilde r^{-k}x)$ converges to $p(\tilde x_0)=x_0$. To see that $\hat i$ is onto the given set, take some sequence $(z_k)_k$ in $X_\infty$ such that $z_k$ converges to $x_0$. Take a neighborhood $V$ of $\tilde x_0$ such that the restriction of $p$ to $V$ is a homeomorphism onto the neighborhood $p(V)$ of $x_0$. Take a smaller neighborhood $U\subset V$ of $\tilde x_0$ such that $\tilde r(U)\subset V$. For $k$ large $z_k$ is in $p(U)$. So $z_k=p(x_k)$ for some $x_k\in U$. Since $r(z_{k+1})=z_k$ it follows that $\tilde r(x_{k+1})=gx_k$ for some $g\in G$. 
But as both $\tilde r(x_{k+1})$ and $x_k$ are in $V$, it follows that $g$ must be $1$. So $x_{k+1}=\tilde r^{-1}x_k$ for $k$ large, bigger than some $k_0$. Then if define $x:=\tilde r^{-k_0}x_{k_0}$,  we have $\hat i(x)=(z_k)_k$.
\par
If $\hat i( y)$ is in $\Omega_{px}$ then $py=px$ so $y=gx$ for some $g\in G$. 
\par
Let us check the $P_x$-measure of the atoms. 
We have 
$$P_{px}(\{\hat i(gx)\})=\lim_{k\rightarrow\infty}P_{px}(\{(z_n)_n\,|\, z_j=p(\tilde r^{-j}x), 0\leq j\leq k\})=
\lim_{k\rightarrow\infty}q^{-k}{m^{(k)}}^*(p\tilde r^{-k}x)h(p\tilde r^{-k}x){m^{(k)}}(p\tilde r^{-k}x)$$
$$=\mathcal P(x)^*h(x_0)\mathcal P(x)=\mathcal P(x)^*\mathcal P(x),$$
and we used the fact from \cite{DuRo1} that the range of $\mathcal P(x)$ is contained in $E_1$ and $h(x_0)=\sum_j h_j(x_0)$ is the projection onto $E_1$. 
\par
Take the cylinder $C_{px,z_1,\dots,z_n}$ of sequences in $\Omega_{px}$ that start with $px,z_1,\dots,z_n$. If we add all the atoms in this cylinder, we obtain, 
\begin{equation}\label{pestea}
\sum_{g, \hat i(gx)\in C_{px,z_1,\dots z_n}}P_{px}(\hat i(gx))=\sum_{g, \hat i(gx)\in C_{px,z_1,\dots z_n}}\mathcal P(x+g)^*\mathcal P(x+g).
\end{equation}
On the other hand , for any section $s\in S$, and with the notation $px_n:=z_n$, for some $x_n\in \tilde X$, using the formula for $h$ in Proposition \ref{duro},
$$\ip{s(z_0)}{P_{px}(C_{px,z_1,\dots,z_n})s(z_0)}=\ip{s(z_0)}{q^{-n}{m^{(n)}}^*(z_n)h(z_n)m^{(n)}(z_n)s(z_0)}$$$$=
\sum_j\sum_{g\in G}q^{-n}\ip{m^{(n)}(z_n)s(z_0)}{\ip{\mathcal P^*(gx_n)s_j(x_0)}{m^{(n)}(z_n)s(z_0)}\mathcal P^*(gx_n)s_j(x_0)}$$
$$=\sum_g\sum_j\left|\ip{\mathcal P^*(gx_n)s_j(x_0)}{q^{-n/2}m^{(n)}(z_n)s(z_0)}\right|^2,$$
and since $\{s_j(x_0)\}_j$ is an orthonormal basis for $E_1$ and the range of $\mathcal P(x)$ is contained in $E_1$,
$$=\sum_g\|\mathcal P(gx_n)q^{-n/2}m^{(n)}(z_n)s(z_0)\|^2.$$
But, from the definition of $\mathcal P$, we have $\mathcal P(gx_n)q^{-n/2}m^{(n)}(z_n)=\mathcal P(g'x)$ for some $g'\in G$ with $\hat i(g'x)\in C_{px,z_1,\dots,z_n}$, and we obtain further
$$=\sum_{g',\hat i(g'x)\in C_{px,z_1,\dots,z_n}}\|\mathcal P(g'x)s(z_0)\|^2.$$
\par
Comparing with \eqref{pestea}, this shows that the sum of the atoms is equal to the measure of the cylinder and thus the measures $P_{px}$ are supported on these atoms. 
\end{proof}

\begin{remark}{\bf Scaling functions.} To obtain old-fashioned scaling functions as described in the introduction, let us consider the case when we are dealing with the Example \ref{mainex}, used in the regular wavelet theory (but the arguments below work also in the more general case we described in this section), and let us consider the case when $\xi$ is the trivial vector bundle $X\times \bc^d=\bt^n\times\bc^d$. 
\par
Then we can take the canonical sections in $\xi$, $c_i(x)=e_i$, $x\in\bt^n$, $i\in\{1,\dots,d\}$, where $e_i$ are the canonical vectors in $\bc^d$. Then, define for each $j\in\{1,\dots,l\}$,
$$\varphi^j_i:=\mathcal W_jc_i,\quad(i\in\{1,\dots,d\}).$$
Then $\varphi^j_i$ is in $\Xi$ so it is a function in $L^2(\br^d)$ (according to \eqref{mutilde}). Also, since $M\mathcal W_j=\W_j$, we have
$$\tilde U\varphi^j_i=\tilde UM\mathcal W_jc_i=\mathcal W_jmc_i\circ\tilde r= \mathcal W_j(\sum_{k=1}^dm_{ki}c_k)=\sum_{k=1}^d\tilde\pi(m_{ki})\varphi^j_k.$$
(Here $\tilde\pi(s)f(x)=s(px)f(x)$, $s\in C(\bt^n), f\in L^2(\br^n), x\in\br^n$).
\par Also, for $f\in C(\bt^n)$,
$$\ip{\varphi^j_i}{\pi(f)\varphi^j_{i'}}=\int_{\br^n}\cj\varphi^j_if\varphi^j_{i'}\,dx=\int_{\bt^n}f\ip{\varphi^j_i}{\varphi^j_{i'}}'\,d\mu=\int_{\bt^n}f\ip{\mathcal W_jc_i}{\mathcal W_jc_{i'}}\,d\mu=\int_{\bt^n}f\ip{c_i}{W_j^*W_jc_{i'}}\,d\mu$$$$=\int_{\bt^n}f\ip{c_i}{h_jc_{i'}}\,d\mu=\int_{\bt^n}f(h_j)_{ii'}d\mu.$$
\par
Thus, $\varphi^j_1,\dots,\varphi^j_n$ form a multi-scaling function in $L^2(\br^n)$, with filter $m$ and correlation matrix $h_j$. 
\par
We can put together all these multi-scaling functions and define
$$\varphi_i:=(\varphi^1_i,\dots,\varphi^l_i)\in\oplus_{j=1}^l L^2(\br^n),\quad(i\in\{1,\dots,d\}).$$
Then we still have the same scaling equation, but now in $\oplus_{j=1}^lL^2(\br^n)$:
$$\tilde U_l\varphi_i=\sum_{k=1}^d\pi_l(m_{ki})\varphi_k,\quad(i\in\{1,\dots,d\}),$$
and the correlation matrix for $\varphi_1,\dots,\varphi_d$ is the harmonic function $h$, i.e., 
$$\ip{\varphi_i}{\pi_l(f)\varphi_{i'}}=\int_{\bt^n}fh_{ii'}\,d\mu,\quad(f\in C(\bt^n)).$$
In Theorem \ref{thj}, we see that if $s=(s^1,\dots,s^d)$, then 
$$Js\circ\theta_0=\oplus_{j=1}^l\mathcal W_js=\oplus_{j=1}^l\mathcal W_j(\sum_{i=1}^ds^ic_i)=\oplus_{j=1}^l(\sum_{i=1}^d\tilde\pi(s^i)\varphi^j_i)=\sum_{i=1}^d\tilde \pi_l(s^i)\varphi_i.$$
This implies, with Theorem \ref{thj}, that the linear span of
$$\{\tilde U_l^j\pi_l(f)\varphi_k\,|\,j\in\bz,f\in C(\bt^n),k\in\{1,\dots,d\}\}$$
is dense in $\oplus_{j=1}^lL^2(\br^n).$
\par
Thus, $\varphi_1,\dots,\varphi_d$ form a ``super''-multi-scaling function for the bigger- (``super-'')space $\oplus_{j=1}^lL^2(\br^n)$.
\end{remark}
\begin{acknowledgements}
We would like to thank professor Palle Jorgensen for his suggestions and for offering us a wider view on the subject, which we included in the introduction.
\end{acknowledgements}
\bibliographystyle{alpha}
\bibliography{dcovariant3}

\end{document}